\documentclass[12pt]{amsart}
\usepackage{amssymb,amsmath,amsthm}
\numberwithin{equation}{section}

\newtheorem{thm}{Theorem}[section]
\newtheorem{lem}[thm]{Lemma}
\newtheorem{cor}[thm]{Corollary}
\newtheorem{ques}[thm]{Question}

\begin{document}
\title[Norms of Weighted Composition Operators]{Reproducing Kernels, de Branges-Rovnyak Spaces, and Norms of Weighted Composition operators}
\author{Michael~T.~Jury} 
\address{Department of Mathematics\\
        University of Florida\\
         Gainesville, Florida 32603}
\email{mjury@math.ufl.edu}
\date{\today}
\begin{abstract} We prove that the norm of a weighted composition operator on the Hardy space $H^2$ of the disk is controlled by the norm of the weight function in the de Branges-Rovnyak space associated to the symbol of the composition operator.  As a corollary we obtain a new proof of the boundedness of composition operators on $H^2$, and recover the standard upper bound for the norm.  Similar arguments apply to weighted Bergman spaces.  We also show that the positivity of a generalized de Branges-Rovnyak kernel is sufficient for the boundedness of a given composition operator on the standard functions spaces on the unit ball.
\end{abstract}
\maketitle

If $H$ is a vector space of functions defined on a set $X$, given a function $b:X\to X$ one can define a \emph{composition operator} $C_b$ by $(C_b f)(x)=f(b(x))$.  When $H$ is the Hardy space $H^2(\mathbb{D})$, the Hilbert space of functions analytic in the open unit disk $\mathbb{D}$ equipped with the norm
$$
\|f\|^2 =\sup_{0< r< 1}\frac{1}{2\pi} \int_0^{2\pi}|f(re^{i\theta})|^2 \, d\theta ,
$$
the composition operator $C_b$ is bounded for every analytic map $b:\mathbb{D}\to \mathbb{D}$, and 
$$
\|C_b\|\leq \left( \frac{1+|b(0)|}{1-|b(0)|}\right)^{1/2}. 
$$
The standard proof of these facts appeals to the Littlewood subordination principle in harmonic analysis; see \cite{MR1397026}.

In this note we give a proof of the boundedness of $C_b$ on $H^2$ which does not use the Littlewood subordination principle, only reproducing kernel methods.  The idea behind the proof is to express the boundedness of certain weighted composition operators $T_f C_b$ in terms of the positivity of kernels related to $b$ and $H^2$, in particular the kernels of the de Branges-Rovnyak spaces.  The boundedness of $C_b$ is obtained as a corollary by a suitable choice of the weight function $f$.  This proof is easily adapted to prove boundedness, with norm estimates, of $C_b$ on the standard weighted Bergman spaces.  We also obtain a sufficient condition for boundedness of composition operators on the standard scale of Hilbert function spaces on the unit ball of $\mathbb{C}^n$.  

We first fix some notation.  For $z,w\in\mathbb{D}$, the \emph{Szeg\H{o} kernel} $k(z,w)$ is defined by
$$
k_w(z)=k(z,w)=\frac{1}{1-\overline{w}z},
$$
$k$ is the reproducing kernel for the Hardy space $H^2$, i.e. for every $f\in H^2$ and every $w\in \mathbb{D}$ we have $\langle f, k_w\rangle =f(w)$. 
Let $b\in H^\infty(\mathbb{D})$ with $\|b\|_\infty \leq 1$, to avoid trivialities we assume $b$ is non-constant.  The \emph{de Branges-Rovnyak space} $H(b)$ is the reproducing kernel Hilbert space on $\mathbb{D}$ with kernel
$$
k^b(z,w)=\frac{1-\overline{b(w)}b(z)}{1-\overline{w}z}.
$$
Equivalently, $H(b)$ may be defined as the linear subspace of $H^2$ equal to the range of the operator ${(I-T_b T_b^*)}^{1/2}$, equipped with the range norm.  The standard references for de Branges-Rovnyak spaces are the books \cite{MR0215065} and \cite{MR1289670}.

For any function $f:\mathbb{D}\to \mathbb{C}$, there is a densely defined operator $T_f^*$ on $H^2$, defined on the Szeg\H{o} kernel $k_w$ by $T_f^*k_w =\overline{f(w)}k_w$ and extended linearly.  (The adjoint notation here is only formal; if $f\in H^\infty$ then $T_f^*$ is bounded and equal to the adjoint of the Toeplitz operator $T_f$.)  We define the operator $C_b^*$ on the linear span of the Szeg\H{o} kernels by $C_b^* k_w =k_{b(w)}$.  If $f\in H^2$ and $f\circ b\in H^2$, then
$$
\langle C_bf, k_w\rangle = \langle f\circ b, k_w \rangle =f(b(w))=\langle f, k_{b(w)}\rangle = \langle f, C_b^* k_w\rangle,
$$
so the operator $C_b^*$ is the formal adjoint of the composition operator $C_b$, so to prove $C_b$ is bounded it suffices to prove that $C_b^*$ is bounded.  It then follows that $C_b^*$ is the genuine adjoint of $C_b$.  

\begin{thm}\label{T:main} For any $f\in H(b)$, the operator $C_b^* T_f^*$ is bounded on $H^2$ and $\|C_b^* T_f^*\|\leq \|f\|_{H(b)}$.
\end{thm}
\begin{proof} We assume $\|f\|_{H(b)} =1$; the general case follows by rescaling.  Put $f_0=f$ and choose unit vectors $f_1, f_2, \dots$ such that $(f_m)_{m\geq 0}$ is an orthonormal basis for $H(b)$.  Then
$$
k^b(z,w)=\frac{1-\overline{b(w)}b(z)}{1-\overline{w}z} =\sum_{m\geq 0} \overline{f_m(w)}f_m(z),
$$
which we rewrite as
\begin{equation}\label{E:main_eqn}
\frac{1}{1-\overline{w}z}= \sum_{m\geq 0} \frac{\overline{f_m(w)}f_m(z)}{1-\overline{b(w)}b(z)}.
\end{equation}
The kernel
\begin{equation}\label{E:main_eqn_2}
\frac{1}{1-\overline{w}z}-\frac{\overline{f(w)}f(z)}{1-\overline{b(w)}b(z)}=\sum_{m\geq 1} \frac{\overline{f_m(w)}f_m(z)}{1-\overline{b(w)}b(z)}
\end{equation}
is therefore positive semidefinite, being a sum of positive semidefinite kernels.  Since we can rewrite the left-hand side as
$$
\langle k_w, k_z\rangle-\langle C_b^* T_f^* k_w, C_b^* T_f^* k_z\rangle
$$
the positivity means that for any $n$ distinct points $w_1, \dots w_n$ in $\mathbb{D}$ and complex numbers $c_1, \dots c_n$, if we define $h\in H^2$ by
$$
h(z)=\sum_{i=1}^n c_i k_{w_i}(z)
$$ 
we then have 
\begin{align}
\nonumber 0&\leq \sum_{i,j=1}^n c_i\overline{c_j} \langle k_{w_i}, k_{w_j}\rangle -\sum_{i,j=1}^n c_i \overline{c_j}\langle C_b^* T_f^* k_{w_i}, C_b^* T_f^* k_{w_j}\rangle \\ \nonumber &=\langle h,h\rangle - \langle C_b^* T_f^* h, C_b^* T_f^* h\rangle,
\end{align}
or $\|C_b^* T_f^* h\|^2\leq \|h\|^2$.  Since such $h$ are dense in $H^2$, it follows that $\|C_b^* T_f^*\|\leq 1$. 
\end{proof}
As a corollary we can now prove that $C_b$ is bounded on $H^2$.  We also note in the proof of Theorem~\ref{T:main} one could deduce the positivity of the left-hand side of equation (\ref{E:main_eqn_2}) directly from the assumption $\|f\|_{H(b)}\leq 1$ (without mention of an orthonormal basis); however the proof given allows us to obtain Corollary~\ref{C:summation} below.
\begin{cor}\label{C:h2_bound} For any analytic map $b:\mathbb{D}\to\mathbb{D}$, the composition operator $C_b$ is bounded on $H^2$, and
$$
\|C_b\|\leq \left( \frac{1+|b(0)|}{1-|b(0)|}\right)^{1/2}.
$$
\end{cor}
\begin{proof}  If $b$ is constant then the boundedness is trivial.  For $b$ non-constant, we apply Theorem 1 to the function
$$
f(z)=k^b_0(z)=1-\overline{b(0)}b(z),
$$
the reproducing kernel for $H(b)$ at the origin.  We have $\|f\|_{H(b)} =(1-|b(0)|^2)^{1/2}$, and we observe that $f$ and $1/f$ are both bounded and analytic in $\mathbb{D}$, since $\|b\|_\infty \leq 1$ and $|b(0)|<1$.  It follows that the Toeplitz operator $T_f^*$ is bounded and invertible with inverse $T_{1/f}^*$, and $\|T_{1/f}^*\|=\|1/f\|_\infty \leq (1-|b(0)|)^{-1}$.  Thus
\begin{align}
\nonumber \|C_b^*\|=\|C_b^* T_f^* T_{1/f}^*\|&\leq \|C_b^* T_f^*\|\|T_{1/f}^*\|\\ \nonumber &\leq \|f\|_{H(b)}\|1/f\|_\infty \\
                                                      \nonumber &\leq \left( \frac{1+|b(0)|}{1-|b(0)|}\right)^{1/2}. 
\end{align}
\end{proof}
\noindent It is well known that the estimate obtained in this corollary is sharp as $b$ ranges over all self-maps $b:\mathbb{D}\to \mathbb{D}$; by a result of Nordgren~\cite{MR0223914} this bound is attained whenever $b$ is an inner function.  Also, it is clear from the above proof that for any $C_b$ we have the norm estimate
\begin{equation}\label{E:hardy_estimate}
\|C_b\| \leq \inf_{f\in H(b)}\left\{ {\|1/f\|}_\infty {\|f\|}_{H(b)}\right\}.
\end{equation}

If we define $T_{f_m}$ in the obvious way on the range of $C_b$, we also obtain immediately from Theorem~\ref{T:main} the following summation identity:
\begin{cor}\label{C:summation}
If $(f_m)_{m\geq 0}$ is any orthonormal basis for $H(b)$, then 
$$
\sum _{m\geq 0} T_{f_m} C_b C_b^* T_{f_m}^* =I
$$
where the sum converges in the strong operator topology.
\end{cor}
\begin{proof}
Equation (\ref{E:main_eqn}) shows that the sum 
\begin{equation}\label{E:op_sum}
\sum _{m\geq 0} T_{f_m} C_b C_b^* T_{f_m}^* 
\end{equation}
converges to $I$ in the weak operator topology.  In fact, the sum converges in the strong operator topology:  subtracting the first $N$ terms of the sum on the right hand side of (\ref{E:main_eqn}) from both sides of that equation shows that the partial sums of (\ref{E:op_sum}) form an increasing sequence of positive operators bounded above (by $I$), so the series is strongly convergent.
\end{proof}
Similar arguments can be used to prove the boundedness of $C_b$ on the standard weighted Bergman spaces $A^2_\alpha$ when $\alpha$ is an integer: for $\alpha\geq 1$ we consider the spaces $A^2_\alpha$ with reproducing kernels 
$$
k^\alpha(z, w)= \frac{1}{(1-z\overline{w})^{\alpha}}.
$$
When $\alpha=1$ this is of course the Szeg\H{o} kernel; for $\alpha>1$ this is the reproducing kernel for the space of analytic functions in $\mathbb{D}$ with norm
$$
\|f\|^2_{A^2_\alpha} =\frac{\alpha-1}{\pi}\int_\mathbb{D}|f(z)|^2 (1-|z|^2)^{\alpha-2} \, dA(z).
$$
For $b$ a self map of the disk as before, we define for integers $\alpha\geq 1$ the space $A(b, \alpha)$ to be the space with reproducing kernel
$$
k^{b, \alpha}(z, w)=\left( \frac{1-b(z)\overline{b(w)}}{1-z\overline{w}}\right)^{\alpha}.
$$
This kernel is positive since it is the $\alpha$-fold Schur product of $k^{b,1}=k^b$ with itself.  Letting $M_f$ denote multiplication by $f$, the same arguments used in the Hardy space prove that if $f\in A(b, \alpha)$ then $M_f C_b$ is bounded on $A^2_\alpha$ and 
\begin{equation}\label{E:alpha_estimate}
\|M_f C_b\|\leq {\|f\|}_{A(b, \alpha)}.
\end{equation}
Applying this inequality to the reproducing kernel for $A(b,\alpha)$ at the origin gives the estimate
$$
{\|C_b\|} \leq \left( \frac{1+|b(0)|}{1-|b(0)|}\right)^{\alpha/2}
$$
and more generally we obtain as in (\ref{E:hardy_estimate})
$$
{\|C_b\|} \leq \inf_{f\in A(b,\alpha)}\left\{ {\|1/f\|}_\infty {\|f\|}_{A(b,\alpha)}\right\}.
$$

It is well known that on the standard Hilbert function spaces on the unit ball $\mathbb{B}^n\subset\mathbb{C}^n$, there are holomorphic maps $b:\mathbb{B}^n\to \mathbb{B}^n$ which do not give bounded composition operators \cite{MR740174, MR783578}.  It is therefore worth understanding why the proof above does not generalize to these spaces.  We consider the spaces $H^2_{n, \alpha}$ for $\alpha\geq 1$, which we define for each $\alpha$ to be the space with reproducing kernel
$$
K^\alpha (z, w)=\frac{1}{(1-\langle z, w\rangle )^\alpha}
$$
where $\langle \cdot, \cdot\rangle$ denotes the standard inner product in $\mathbb{C}^n$.  When $\alpha =n$ this is the Hardy space on $\mathbb{B}^n$, and when $\alpha=n+1$ this is the Bergman space.  If we attempt to adapt the single-variable argument to this setting, we are led to consider in place of the de Branges-Rovnyak kernel the kernel
$$
K^{b,\alpha}(z, w)=\left( \frac{1-\langle b(z), b(w)\rangle}{1-\langle z, w\rangle}\right)^\alpha.
$$
However, even when $\alpha$ is an integer, this kernel is not positive semidefinite for all holomorphic maps $b:\mathbb{B}^n\to \mathbb{B}^n$.  (When $\alpha=1$, this occurs because there exist holomorphic functions on $\mathbb{B}^n$ bounded by $1$ that do not act as contractive multipliers of the the Hilbert function space $H^2_{n,1}$; see e.g.\,\cite[Chapter 8]{MR1882259}.)  Nonetheless, whenever this kernel is positive, the composition operator $C_b$ is bounded on $H^2_{n,\alpha}$ and we obtain a norm estimate analogous to the one-variable case.  This could be proved by a modification of the argument in the one variable case, except for one complication:  in the proof of Corollary 2, we used the fact the the reciprocal of the reproducing kernel for $H(b)$ induces a bounded multiplication operator on $H^2$; this follows simply because the reciprocal of the kernel is bounded holomorphic function.  However, on the ball (when $\alpha <n$) the boundedness in the supremum norm is not sufficient to give a bounded multiplication operator, so an extra argument, provided by the following lemma, is necessary.  For a holomorphic map $b:\mathbb{B}^n\to \mathbb{B}^n$, we write $b=(b_1, \dots b_n)$ for the coordinate functions of $b$, so each $b_i$ is a holomorphic map from $\mathbb{B}^n$ to the unit disk and $\sum_{i=1}^n |b(z)|^2 <1$ for all $z\in\mathbb{B}^n$.  Finally, for any $c=(c_1, \dots c_n)\in \mathbb{C}^n$ we write $|c|=(\sum_{i=1}^n |c_i|^2)^{1/2}$.
\begin{lem}\label{L:one_over_k}  Let $b=(b_1\dots b_n)$ be a holomorphic map from $\mathbb{B}^n$ into itself and let $\alpha\geq 1$.  If the kernel
$$
K^{b,\alpha}(z, w)=\left( \frac{1-\langle b(z), b(w)\rangle}{1-\langle z, w\rangle}\right)^\alpha
$$
is positive semidefinite on $\mathbb{B}^n$, then each coordinate function $b_i$ is a contractive multiplier of $H^2_{n, \alpha}$.  Moreover the function 
$$
\frac{1}{K^{b,\alpha}(z,0)}=(1-\langle b(z) , b(0)\rangle)^{-\alpha}
$$
is a bounded multiplier of $H^2_{n,\alpha}$ of norm at most $(1-|b(0)|)^{-\alpha}$.
\end{lem}
\begin{proof} Fix $\alpha\geq 1$ and let $M_{b_i}$ denote the operator of multiplication by $b_i$ on $H^2_{n,\alpha}$.  Assuming $K^{b,\alpha}$ is positive, the kernel
$$
\frac{1-\langle b(z), b(w)\rangle}{(1-\langle z, w\rangle)^\alpha}=\frac{1}{(1-\langle b(z),b(w)\rangle)^{\alpha-1}} K^{b,\alpha}(z,w)
$$
is positive semidefinite, since it is the Schur product of the positive kernels $K^{b,\alpha}$ and $(1-\langle b(z),b(w)\rangle)^{1-\alpha}$.  By a standard reproducing kernel argument, this kernel is positive if and only if the operator $I-\sum_{i=1}^n M_{b_i}M_{b_i}^*$ is positive.  Thus each $M_{b_i}$ is contractive.  

To prove the second statement, we first observe that for each $w\in\mathbb{B}^n$, 
$$
\|\sum_{i=1}^n \overline{b_i(w)}M_{b_i}\|\leq |b(w)|.
$$
To see this, note that the operator inequality
$$
I-\sum_{i=1}^n M_{b_i}M_{b_i}^* \geq 0
$$
may be interpreted as saying the column operator $(M_{b_1} \cdots M_{b_n})^*$ is contractive from $H^2_{n,\alpha}$ to the direct sum of $n$ copies of this space with itself; and similarly for the row operator $(M_{b_1} \cdots M_{b_n})$ in the reverse direction.  Since for any $c=(c_1, \dots c_n)\in\mathbb{C}^n$, the column operator $(c_1 I \cdots c_nI)^T$ has norm $|c|$, we have 
$$
\|\sum_{i=1}^n \overline{b_i(w)}M_{b_i}\|=\|(M_{b_1} \cdots M_{b_n})(b_1(w) I \cdots b_n(w)I)^*\|\leq |b(w)|.
$$ 
The second claim of the lemma can now be proved by expanding the operator $(I-\sum_{i=1}^n \overline{b_i(0)}M_{b_i})^{-\alpha}$ as a power series in $\sum_{i=1}^n \overline{b_i(0)}M_{b_i}$, which is norm convergent since $\|\sum_{i=1}^n \overline{b_i(0)}M_{b_i}\|\leq |b(0)|<1$.
\end{proof}

\begin{thm}\label{T:ball} Let $b:\mathbb{B}^n\to \mathbb{B}^n$ be a holomorphic map and let $\alpha\geq 1$. Suppose the kernel $K^{b,\alpha}(z,w)$ is positive semidefinite on $\mathbb{B}^n$, and let $H^2_{n,\alpha}(b)$ denote the Hilbert space with kernel $K^{b,\alpha}$.   Then for all $f\in H^2_{n,\alpha}(b)$, the weighted composition operator $M_f C_b$ is bounded on $H^2_{n,\alpha}$ and 
$$
\|M_f C_b \|\leq \|f\|_{H^2_{n,\alpha}(b)}.
$$
Moreover $C_b$ is bounded on $H^2_{n,\alpha}$ and 
$$
\|C_b\| \leq \left(\frac{1+|b(0)|}{1-|b(0)|}\right)^{\alpha/2}.
$$
\end{thm}
\begin{proof} The first inequality is proved exactly as in the one-variable case.  Similarly, by virtue of Lemma~\ref{L:one_over_k} the norm estimate for $C_b$ follows from the estimate for the weighted operator applied to $f(z)=K^{b,\alpha}(z,0) $ as in Corollary~\ref{C:h2_bound}.
\end{proof}
For certain values of $\alpha$ (e.g $\alpha=n$, the Hardy space) there are known necessary and sufficient conditions for the boundedness of $C_b$ on $H^2_{n,\alpha}$, given in terms of Carleson measures \cite{MR783578}.  Theorem~\ref{T:ball} tells us that positivity of $K^{b,\alpha}$ is sufficient for the boundedness of $C_b$ on $H^2_{n,\alpha}$.  In the one-variable Hardy space, this condition is also necessary:  the symbol of a bounded composition operator on $H^2$ must belong to the unit ball of $H^\infty$, which coincides with the unit ball of the multiplier algebra of $H^2$; the positivity of the de~Branges-Rovnyak kernel follows.  However in general the positivity of $K^{b,\alpha}$ is not necessary for the boundedness of $C_b$.  For example, when $n=2$ and $\alpha=1$ it can be shown by standard estimates that for 
\begin{equation*}
b_r(z_1, z_2)=(2rz_1 z_2, 0) 
\end{equation*}
the composition operator $C_{b_r}$ is bounded when $0\leq r<1$ and unbounded when $r=1$.  If the kernels $K^{b_r, 1}$ were positive for all $r<1$ then by taking pointwise limits $K^{b_1, 1}$ would be positive as well, which by Theorem~\ref{T:ball} would make $C_{b_1}$ bounded, a contradiction.  Thus $K^{b_r,1}$ is non-positive for $r$ sufficiently close to $1$.

In the ball, when $\alpha=1$, the positivity of $K^{b,1}$ is equivalent to the statement that the tuple $(M_{b_1}, \dots M_{b_n})$ is a row contraction.  Recently, S. Shimorin~\cite{MR2139107} has proved that the positivity of $K^{b,1}$ is essentially equivalent to a commutant lifting theorem between $H^2_{n,1}$ and the space with reproducing kernel $(1-\langle b(z),b(w)\rangle)^{-1}$.  Also, it can be shown that the linear fractional maps of the unit ball introduced in \cite{MR1768872} induce row contractions; we thus obtain a new proof of the boundedness of linear fractional composition operators on the ball, with norm estimates.  This will be discussed in detail in a separate paper.

Finally, returning to Theorem 1, it is clear $f\in H(b)$ is only a sufficient condition for the boundedness of $T_f C_b$ on $H^2$, not a necessary one; e.g. $T_fC_b$ is bounded for all $f\in H^\infty$ but in general $H(b)$ does not contain $H^\infty$.  It follows from the theorem that $T_f C_b$ is bounded for all $f$ in the linear span of $H^\infty \cdot H(b)$; however it is not clear how close this comes to describing all bounded weighted composition operators (with analytic weights).  In particular, we do not know if this set of weights is ``dense'' in the set of all weights giving bounded operators, in the following sense:
\begin{ques}  Given $f\in H^2$ and $b$ in the unit ball of $H^\infty$ such that $T_f C_b$ is bounded on $H^2$, and given $\epsilon >0$, does there exists $g\in \text{span}\,(H^\infty \cdot H(b))$ such that $\|T_{f-g} C_b\|<\epsilon$?
\end{ques}
A necessary and sufficient condition for the boundedness of $T_fC_b$ can be given in terms of Carleson measures \cite[Theorem 2.2]{MR1864316}, but the relationship between such conditions and the $H(b)$ spaces is not immediately clear.

\bibliographystyle{plain} 
\bibliography{comp_dBR} 
\end{document}